\documentclass[12pt]{article}

\usepackage{amsmath, epsfig, cite, lineno}
\usepackage{amssymb,amsthm}
\usepackage{amsfonts, color}
\usepackage{latexsym}

\newtheorem{thm}{Theorem}[section]

\newtheorem{conj}[thm]{Conjecture}
\newtheorem{lem}[thm]{Lemma}

\numberwithin{equation}{section}

\begin{document}

%\linenumbers

\begin{center}
{\Large\bf On Zudilin's $q$-question about Schmidt's problem}
\end{center}

\vskip 2mm \centerline{Victor J. W. Guo$^1$  and Jiang Zeng$^{2}$}
\begin{center}
{\footnotesize $^1$Department of Mathematics, East China Normal University,\\ Shanghai 200062,
 People's Republic of China\\
{\tt jwguo@math.ecnu.edu.cn,\quad http://math.ecnu.edu.cn/\textasciitilde{jwguo}}\\[10pt]
%$^2$Universit\'e de Lyon, Lyon, F-69003, France\\
$^2$Universit\'e de Lyon; Universit\'e Lyon 1; Institut Camille
Jordan, UMR 5208 du CNRS;\\ 43, boulevard du 11 novembre 1918,
F-69622 Villeurbanne Cedex, France\\
{\tt zeng@math.univ-lyon1.fr,\quad
http://math.univ-lyon1.fr/\textasciitilde{zeng}} }
\end{center}

%%date: January 4, 2011
%\vskip 5mm
%\noindent {\it Suggested Running title}: Two Identities of Gould

\vskip 0.7cm \noindent{\bf Abstract.}
We propose an elemantary approach to Zudilin's $q$-question about Schmidt's problem
[Electron. J. Combin. 11 (2004), \#R22], which has been solved in
a previous paper [Acta Arith. 127 (2007), 17--31].  The new approach is
based on a  $q$-analogue of  our recent result  in  [J. Number Theory 132 (2012), 1731--1740]
derived from $q$-Pfaff-Saalsch\"utz identity.



\section{Introduction}

%%%%%%%%%%%%%%%%%%%%%%%%%%%%%%%%%%%%%%%%%%%%%%%%%%%%%%%%%%%%%%%%%%%%%%%%%%%%%%%%%%%%%%%%%%%%%%%%%%%%%

In 2007,  answering a question of Zudilin \cite{Zudilin}, the following result was proved in \cite{GJZ}.
\begin{thm}\label{thm:GJZ}
Let $r\geqslant  1$. Then there exists a unique sequence of polynomials $\{c_{i}^{(r)}(q)\}_{i=0}^\infty$
in $q$ with nonnegative integral coefficients such that, for any $n\geq 0$,
\begin{align}\label{qzudilin}
&\sum_{k=0}^nq^{r{n-k\choose 2}+(1-r){n\choose 2}}
{n\brack k}^r{n+k\brack k}^r
=\sum_{i=0}^nq^{{n-i\choose 2}+(1-r){i\choose 2}}
{n\brack i}{n+i\brack i}c^{(r)}_i(q).
\end{align}
\end{thm}
Here,  the $q$-binomial coefficients ${n\brack k}$  are defined by
$$
{n\brack k}
=\begin{cases}\displaystyle\frac{(q)_n}{(q)_k (q)_{n-k}}, &\text{if $0\leqslant  k\leqslant  n$}, \\[5pt]
0, &\text{otherwise,}
\end{cases}
$$
where $(q)_0=1$ and $(q)_n=(1-q)(1-q^2)\cdots (1-q^n)$ for $n=1,2,\ldots.$ It is well known that
${n\brack k}$ is a polynomial in $q$ with nonnegative integral coefficients of degree $k(n-k)$
(see \cite[p.~33]{Andrews98}).

 The proof of
\eqref{qzudilin} given in \cite{GJZ} is a $q$-analogue of Zudilin's \cite{Zudilin}
approach to Schmidt's problem  (see \cite{Schmidt, Strehl})  by first using  the $q$-Legendre inversion formula to obtain a formula for $c^{(r)}_k(q)$ and then
applying  a basic hypergeometric identity due to Andrews~\cite{Andrews75} to show
that the latter expression is indeed a polynomial in $q$ with nonnegative integral coefficients.
In this paper we propose a new and elementary approach to Zudilin's $q$-question,
which yields
not only  a new proof of Theorem \ref{thm:GJZ}, but also
more solutions to Zudilin's $q$-question about  Schmidt's problem.

 Our  starting point is  the following $q$-version of Lemma 4.2 in \cite{GZ}.
\begin{lem}\label{thm:pkir}
Let $k\geqslant 0$ and $r\geqslant 1$. Then there exists a unique sequence of Laurent polynomials
$\{P_{k,i}^{(r)}(q)\}_{i=k}^{rk}$
in $q$ with nonnegative integral coefficients such that, for any $n\geqslant k$,
\begin{align}
{n\brack k}^r{n+k\brack k}^r
=\sum_{i=k}^{\min\{n,rk\}} q^{(rk-i)n}{n\brack i}{n+i\brack i}P_{k,i}^{(r)}(q). \label{eq:q-nkrnk}
\end{align}
Moreover, the polynomials $P_{k,i}^{(r)}(q)$ can be computed recursively  by $P_{k,k}^{(1)}(q)=1$ and
\begin{align}
P_{k,k+j}^{(r+1)}(q)=\sum_{i=k}^{rk}q^{(j-i)(j+k)}
{k+i\brack i}{k\brack i-j}{k+j\brack j}P_{k,i}^{(r)}(q),\ 0\leqslant  j\leqslant  rk.  \label{eq:pkkj}
\end{align}
\end{lem}

To derive Theorem 1.1 from Lemma 1.2 we first consider a more general problem.
Let $f(x,y)$ and  $g(x,y)$ be any  polynomials in $x$ and $y$ with integral coefficients.
Multiplying  \eqref{eq:q-nkrnk} by $q^{-nkr+f(k,r)}$ and summing  over $k$ from $0$ to $n$
we obtain
\begin{align}
\sum_{k=0}^n q^{-nkr+f(k,r)} {n\brack k}^r {n+k\brack k}^r
=\sum_{i=0}^n q^{-ni-g(i,r)}
{n\brack i}{n+i\brack i}\sum_{k=0}^iT_{k,i}^{(r)}(q), \label{eq:qProb2master}
\end{align}
where
%$S_{i}^{(r)}(q)=\sum_{k=0}^iT_{k,i}^{(r)}(q)$with
\begin{align}\label{eq:T}
T_{k,i}^{(r)}(q)=q^{f(k,r)+g(i,r)}P_{k,i}^{(r)}(q),\ 0\leqslant k\leqslant i,
\text{ and $P_{k,i}^{(r)}(q)=0$ if $i>kr$.}
\end{align}
Obviously $T_{k,i}^{(r)}(q)$ are Laurent polynomials in $q$ with nonnegative integral coefficients.
For example,  taking
$f=g=0$, we immediately obtain the following result.

\begin{thm}\label{thm:qProb2}
Let $r\geqslant  1$. Then there exists a unique  sequence of Laurent
polynomials $\{b_{i}^{(r)}(q)\}_{i=0}^{\infty}$ in $q$ with nonnegative integral coefficients such that,
for any $n\geq 0$,
\begin{align}
\sum_{k=0}^n q^{-rkn} {n\brack k}^r {n+k\brack k}^r
=\sum_{i=0}^n q^{-ni}
{n\brack i}{n+i\brack i}b_{i}^{(r)}(q). \label{eq:qProb2}
\end{align}
Moreover, we have $b^{(r)}_i(q)=\sum_{k=0}^i P_{k,i}^{(r)}(q)$.
\end{thm}

Now, we look for a sufficient condition for $T_{k,i}^{(r)}(q)$ in
\eqref{eq:qProb2master} to be a polynomial. It follows from
\eqref{eq:pkkj}  that
\begin{align}\label{eq:stepT}
T_{k,i}^{(r+1)}(q)=  \sum_{j=k}^{rk}q^{A} {k+j\brack j}{k\brack i-j}  {i\brack k}T_{k,j}^{(r)}(q),
\end{align}
where
\begin{align}\label{eq:gen}
A=f(k,r+1)+g(i,r+1)-f(k,r)-g(j,r)+i(i-k-j).
\end{align}
Hence,  the positivity of $A$ will ensure that  $T_{k,i}^{(r)}(q)$  is  a polynomial in $q$.

 We shall first prove Lemma~\ref{thm:pkir} in the next section and
 then prove Theorem~\ref{thm:GJZ} in Section 3 by choosing special polynomials $f$ and $g$. Some open problems are raised in Section 4.

 %%%%%%%%%%%%%%%%%%%%%%
\section{Proof of Lemma~\ref{thm:pkir}}
 We proceed by induction on $r$.
We need the following form of  Jackson's
$q$-Pfaff-Saalsch\"utz identity (see \cite[pp. 37-38]{Andrews98} or \cite{Schmidt}  for example):
\begin{align}
{m+n\brack M}{n\brack N}=\sum_{j\geqslant 0} q^{(N-j)(M-m-j)}
{M-m\brack j}{N+m\brack m+j}{m+n+j\brack M+N}. \label{eq:qPfaff}
\end{align}
Substituting $m\to k-i$, $n\to n+i$, $M\to n-i$ and $N\to i$ in \eqref{eq:qPfaff},
we get
\begin{align*}
{n+k\brack n-i}{n+i\brack i}=\sum_{j=0}^{i}q^{(i-j)(n-k-j)}{n-k\brack j}{k\brack i-j}{n+k+j\brack n},
\end{align*}
which can be rewritten as
\begin{align}
{n\brack i}{n+i\brack i}
=\sum_{i=0}^i q^{(i-j)(n-k-j)}\frac{(q)_{k+i} (q)_{j}}{(q)_{k+j}(q)_{i}}{k\brack i-j}{n-k\brack j}{n+k+j\brack j}.
\label{eq:lem03}
\end{align}

 It is clear that \eqref{eq:q-nkrnk} holds for $r=1$ with
$P_{k,k}^{(r)}(q)=1$.
Suppose that  \eqref{eq:q-nkrnk} holds for some $r\geqslant  1$.
Multiplying both sides of \eqref{eq:q-nkrnk} by
${n\brack k} {n+k\brack k}$ and applying \eqref{eq:lem03}, we immediately get
\begin{align}
{n\brack k}^{r+1} {n+k\brack k}^{r+1}
&=\sum_{i=k}^{rk} q^{(rk-i)n} {n\brack k}{n+k\brack k}P_{k,i}^{(r)}(q) \nonumber \\
&\hskip 7mm \times \sum_{j=0}^{i}  q^{(i-j)(n-k-j)}\frac{(q)_{k+i} (q)_{j}}{(q)_{k+j}(q)_{i}}
{k\brack i-j}{n-k\brack j}{n+k+j\brack j}
 \nonumber \\
&=\sum_{j=0}^{rk}  q^{(rk-j)n}{n\brack k+j}{n+k+j\brack k+j}
P_{k,k+j}^{(r+1)}(q),  \label{eq:pkir}
\end{align}
where $P_{k,k+j}^{(r+1)}(q)$ is given by \eqref{eq:pkkj}.
By the induction hypothesis, these $P_{k,k+j}^{(r+1)}(q)$ are Laurent polynomials in $q$ with
nonnegative integral coefficients.
Hence Lemma \ref{thm:pkir} is true for $r+1$.
%%%%%%%%%%%%%%%
\section{Proof of Theorem~\ref{thm:GJZ}}
In \eqref{eq:qProb2master}, taking
$ f(k,r)=r{k+1\choose 2}$,  $g(i,r)=(r-2){i\choose 2}-i$,
and multiplying by $q^{n\choose 2}$, we obtain \eqref{qzudilin} with
\begin{align}
c^{(r)}_i(q)=q^{(r-2){i\choose 2}-i}\sum_{k=0}^i q^{r{k+1\choose 2}}P_{k,i}^{(r)}(q). \label{eq:brkq}
\end{align}
 By  \eqref{eq:gen}  the corresponding $A$ reads as follows
$$
A=(r-2)\left[{i\choose 2}-{j\choose 2}\right]+{i-k\choose 2}+(i-1)(i-j).
$$
If $r\geqslant 2$, since $i\geqslant j$, we have $A\geqslant 0$.
If $r=1$, then \eqref{eq:stepT} implies that $j=k$ and
$A=2{i-k\choose 2}\geqslant 0$.
Thus the
$c^{(r)}_i(q)$ in \eqref{eq:brkq}  is a polynomial in $q$.
For example,   by \eqref{eq:T} we have
$$
T_{k,i}^{(2)}(q)=  q^{2{i-k\choose 2}}{2k\brack i}{i\brack k}^2,
$$
and
$$
c^{(2)}_i(q)=\sum_{k=0}^i q^{2{i-k\choose 2}}{2k\brack i}{i\brack k}^2,
$$
which coincides with \cite[(3,1)]{GJZ}.

\section{Open problems}
For any positive integers $r$ and $s$,  it is easy to see that there are uniquely determined
rational numbers  $c_k^{(r,s)}$ ($k\geqslant 0$), independent of $n$ ($n\geqslant 0$),  satisfying
\begin{equation}\label{eq:zu}
\sum_{k=0}^n{n\choose k}^r{n+k\choose k}^r =\sum_{k=0}^n{n\choose k}^s{n+k\choose k}^s c_k^{(r,s)}.
\end{equation}
When $s=1$ and $r\geqslant 1$,  the integrality of $c_k^{(r,s)}$ is the original problem of Schmidt~\cite{Schmidt}.
When $s>1$ and $r>s$, we observe that the numbers $c_k^{(r,s)}$  are not always integers.  From  arithmetical point of view, the following problems
may be interesting.

\begin{conj}\label{prob:new2}
For any $s>1$ and $n\geqslant 0$,  there is an integer $r>s$ such that all the numbers $c_{k}^{(r,s)}\ (0\leqslant k\leqslant n)$
are integers.
\end{conj}

For $s=2$,  via Maple, we find  that  the least such integers $r:=r(n,s)$ are
$r(0,2)=r(1,2)=r(2,2)=3,r(3,2)=7,r(4,2)=32,r(5,2)=212$.

\begin{conj}\label{prob:new1}
For any $r>s>1$, there is a positive integer $n$ such that $c_{n}^{(r,s)}$
is not an integer.
\end{conj}

\vskip 5mm \noindent{\bf Acknowledgments.} This work was partially supported by the Fundamental
Research Funds for the Central Universities.

\end{document}